\documentclass[11pt]{article}
\usepackage{mathrsfs}
\usepackage{amsthm}
\usepackage{amssymb}
\usepackage{amsmath}
\usepackage{graphicx}
\usepackage{color}
\usepackage{amsfonts}
\usepackage{float}
\usepackage{cite}
\usepackage[text={140mm,210mm},left=35mm,vmarginratio=1:1]{geometry}
\newtheorem{theorem}{Theorem}[section]
\newtheorem{remark}{Remark}[section]

\newtheorem{proposition}[theorem]{Proposition}

\numberwithin{equation}{section}
\normalsize

\begin{document}
\title{\textbf{Statistical inference for unknown parameters of stochastic SIS epidemics on complete graphs}}
\author{Huazheng Bu \thanks{\textbf{E-mail}: 19121617@bjtu.edu.cn \textbf{Address}: School of Science, Beijing Jiaotong University, Beijing 100044, China.} and { Xiaofeng Xue} \thanks{\textbf{E-mail}: xfxue@bjtu.edu.cn \textbf{Address}: School of Science, Beijing Jiaotong University, Beijing 100044, China.}\\ Beijing Jiaotong University}
\date{}
\maketitle

\noindent {\bf Abstract:}
In this paper, we are concerned with the stochastic susceptible-infectious-susceptible (SIS) epidemic model on the complete graph with $n$ vertices. This model has two parameters, which are the infection rate and the recovery rate. By utilizing the theory of density-dependent Markov chains, we give consistent estimations of the above two parameters as $n$ grows to infinity according to the sample path of the model in a finite time interval. Furthermore, we establish the central limit theorem (CLT) and the moderate deviation principle (MDP) of our estimations. As an application of our CLT, reject regions of hypothesis testings of two parameters are given. As an application of our MDP, confidence intervals with lengths converging to $0$ while confidence levels converging to $1$ are given as $n$ grows to infinity.

\quad


\section{Introduction}\label{section one}

In this paper, we are concerned with statistical inference for unknown parameters of stochastic susceptible-infectious-susceptible (SIS) epidemics on complete graphs with large degrees. For integer $n\geq 2$, a complete graph $C_n$ is a graph with $n$ vertices, where any two vertices are connected with an edge. A stochastic SIS epidemic model on $C_n$, which is also named as a contact process (see Section 6 of Liggett's book \cite{Lig1985}), is a continuous time Markov process with state space $\wp(C_n):=\{A:~A\subseteq C_n\}$. Let $A_t$ be the state of the process at moment $t$ for $t\geq 0$, then $\{A_t\}_{t\geq 0}$ evolves as follows. For any $x\in C_n$ and $t\geq 0$,
\[
A_t\text{~flips to~}
\begin{cases}
A_t\setminus \{x\} \text{~at rate~}\theta & \text{~if~}x\in A_t,\\
A_t\cup \{x\} \text{~at rate~} \frac{\lambda}{n}|A_t| & \text{~if~}x\not \in A_t,
\end{cases}
\]
where $|A|$ is the cardinality of a set $A$ while $\lambda, \theta$ are two parameters called `infection rate' and `recovery rate' respectively. Note that an event occurs at rate $r$ for some $r>0$ means that the random time we wait for the event to occur follows exponential distribution with parameter $r$.

Intuitively, $\{A_t\}_{t\geq 0}$ describes the spread of a susceptible-infectious-susceptible epidemic on $C_n$. Vertices in $A_t$ are infectious at moment $t$ while those out of $A_t$ are susceptible. An infectious vertex becomes susceptible at rate $\theta$ while a susceptible vertex is infected at rate proportional to the number of infectious neighbours, which is $|A_t|$ since any two vertices on the complete graph are neighbours.

In this paper, we consider $\theta$ and $\lambda$ as unknown parameters which do not rely on $n$. The aim of this paper is to give consistent estimations $\widehat{\lambda}$ and $\widehat{\theta}$ of $\lambda$ and $\theta$ respectively by observing the path of $\{A_t\}_{0\leq t\leq T_0}$ for large $n$ and given moment $T_0>0$. Furthermore, we will establish central limit theorem and moderate deviation principle for $\widehat{\lambda}$ and $\widehat{\theta}$. For mathematical results and their applications, see next section.

Note that $C_n$ is a finite graph while $\emptyset$ is an absorbed state of the process. If we fix $n$ while let $t$ grow to infinity, we can only find that all the vertices are susceptible eventually. That's why in our setting we fix the moment $T_0$ while let the scale of the graph grow to infinity. The other setting where the graph is fixed while time $t$ grows to infinity can be investigated when the graph is infinite. For related literatures, see References \cite{Becker1977, Fierro2015, Guy2015, Hadeler2011, Lekone2006, Lin2013, Pan2014, Yip1998} and so on.

\section{Main results and their applications}\label{section two}

In this section we give our main results and some of their applications. From now on we let $T_0>0$ be a fixed given moment and assume that
\[
A_0=C_n,
\]
i.e., all the vertices are infectious initially. First we give consistent estimations of $\lambda$ and $\theta$. For this purpose, we define
\[
K(x,y)=
\begin{cases}
\frac{1}{1+xT_0} & \text{~if~}x=y,\\
\frac{(y-x)e^{(y-x)T_0}}{ye^{(y-x)T_0}-x} & \text{~else}
\end{cases}
\]
for $x, y>0$. It is easy to check that $K(x,y)$ is continuous and strictly increasing with $y$, so it is reasonable to define $H(x,z)=H(x, \cdot)(z)$ as the inverse function of $K(x, \cdot)(y)$, i.e.,
\[
H(x,z)=y \text{~if and only if~}K(x,y)=z.
\]
When we need to distinguish different $C_n$s, we write $A_t$ as $A_t^n$. For each $n\geq 1$, we define
\[
X_n=\left|A_{T_0}^n\right| \text{~and~}V_n=\left|\left\{x\in C_n:~x\in A_t^n \text{~for all~}0\leq t\leq T_0\right\}\right|,
\]
where $|A|$ is the cardinality of the set $A$. That is to say, $X_n$ is the number of infectious vertices at moment $T_0$ while $V_n$ is the number of vertices maintaining infectious during $[0, T_0]$. Note that $X_n$ and $V_n$ are statistics which can be observed directly according to the trajectory of $\{A_t^n\}_{0\leq t\leq T_0}$.

We define
\[
\begin{cases}
\hat{\theta}_n=-\frac{1}{T_0}\log{\frac{V_n}{n}},\\
\hat{\lambda}_n=H\left(\hat{\theta}_n, \frac{X_n}{n}\right),
\end{cases}
\]
then we have the following result, which gives consistent estimations of $\theta$ and $\lambda$.

\begin{theorem}\label{theorem 2.1 LLN}
Under the assumption that $A_0^n=C_n$ for all $n\geq 1$,
\[
\lim_{n\rightarrow}\hat{\theta}_n=\theta \text{~and~} \lim_{n\rightarrow+\infty}\hat{\lambda}_n=\lambda
\]
in probability.
\end{theorem}
By Theorem \ref{theorem 2.1 LLN}, for the contact process on $C_n$ where $n$ is large and all the vertices are infectious initially, we can give estimations of $\lambda$ and $\theta$ with small errors by observing the trajectory of $\{A_t^n\}_{0\leq t\leq T_0}$ and then recording $X_n$ and $V_n$. Note that the advantage of this approach is that we do not need to observe this contact process for a long time. Following are simulation results of $\hat{\lambda}_n, \hat{\theta}_n, $ under four different settings of $\lambda, \theta$ for $20\leq n \leq 1000$ and $T_0=1$.

\textbf{Setting 1} $\lambda=0.3$ and $\theta=1$.
\begin{figure}[H]
\centering
\includegraphics[scale=0.5]{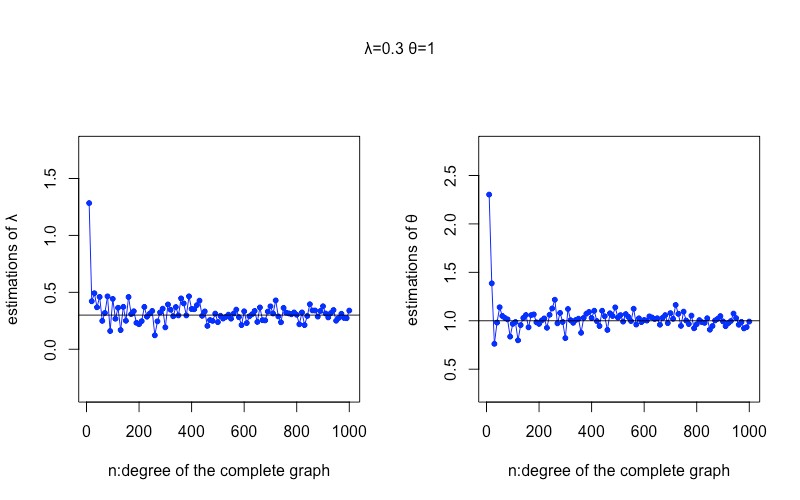}
\caption{$\lambda=0.3$ and $\theta=1$}\label{figure 2.1}
\end{figure}

\textbf{Setting 2} $\lambda=0.5$ and $\theta=1$.
\begin{figure}[H]
\centering
\includegraphics[scale=0.5]{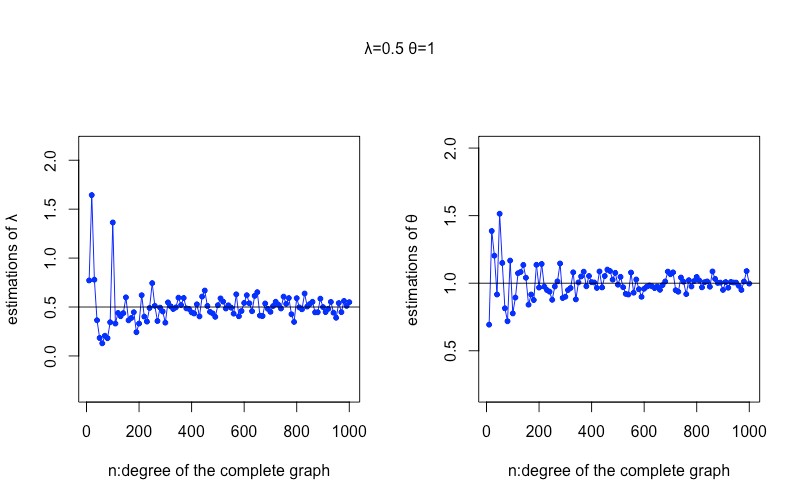}
\caption{$\lambda=0.5$ and $\theta=1$}\label{figure 2.2}
\end{figure}

\textbf{Setting 3} $\lambda=1$ and $\theta=1$.
\begin{figure}[H]
\centering
\includegraphics[scale=0.5]{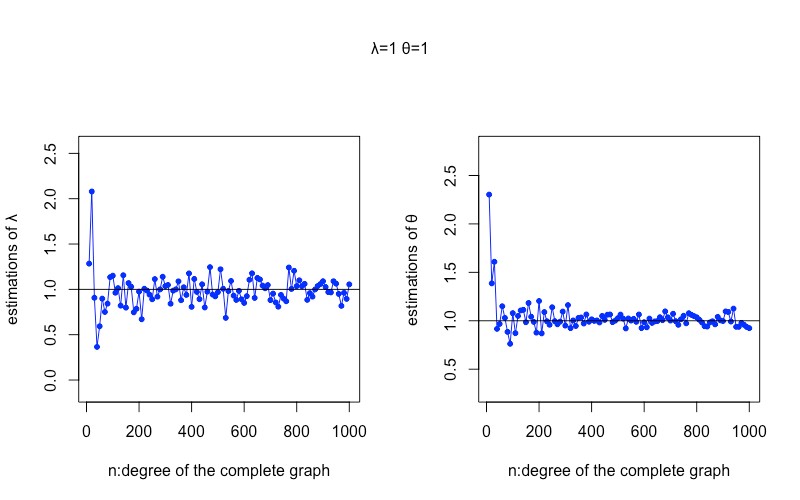}
\caption{$\lambda=1$ and $\theta=1$}\label{figure 2.3}
\end{figure}

\textbf{Setting 4} $\lambda=2$ and $\theta=1$.
\begin{figure}[H]
\centering
\includegraphics[scale=0.5]{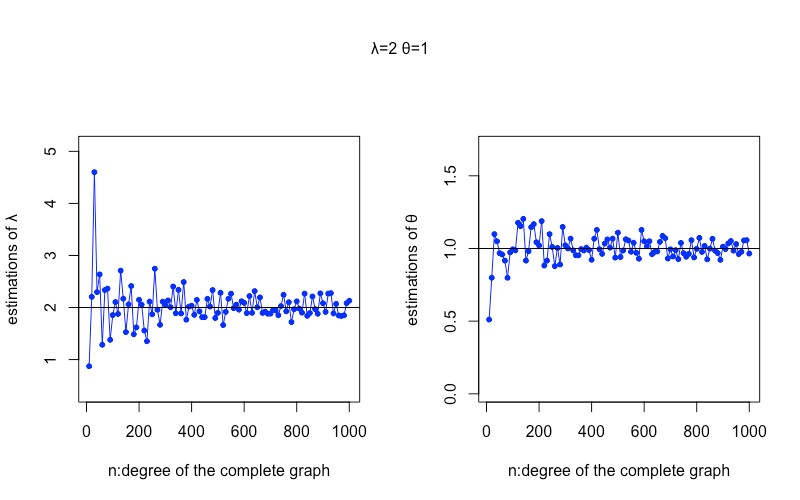}
\caption{$\lambda=2$ and $\theta=1$}\label{figure 2.4}
\end{figure}

Our next theorem gives the central limit theorem of $\hat{\lambda}_n$ and $\hat{\theta}_n$. To give our result, we need to introduce some notations and definitions.
Let $M_1, M_2$ be $2\times 2$ matrices defined as
\[
M_1=\begin{pmatrix}
\frac{\partial H}{\partial z}\left(\theta, K(\theta, \lambda)\right) & \frac{\partial H}{\partial x}\left(\theta, K(\theta, \lambda)\right)\\
0 & 1
\end{pmatrix}
\text{~and~}
M_2=\begin{pmatrix}
1 & 0 \\
0 & -\frac{e^{\theta T_0}}{T_0}
\end{pmatrix}
.
\]
We then use $M_3$ to denote $M_1M_2$, i.e.,
\[
M_3=
\begin{pmatrix}
\frac{\partial H}{\partial z}\left(\theta, K(\theta, \lambda)\right) & -\frac{\partial H}{\partial x}\left(\theta, K(\theta, \lambda)\right)\frac{e^{\theta T_0}}{T_0}\\
0 &  -\frac{e^{\theta T_0}}{T_0}
\end{pmatrix}
.
\]
We use $\mathsf{T}$ to denote the transposition operator. For later use, we define
\[
l_1=(-1, -1)^{\mathsf{T}}, l_2=(-1, 0)^{\mathsf{T}} \text{~and~}l_3=(1, 0)^{\mathsf{T}},
\]
while
\[
F_1(x,y)=\theta y, F_2(x,y)=\theta(x-y) \text{~and~} F_3(x,y)=\lambda x(1-x).
\]
For $0\leq t\leq T_0$ and our two unknown parameters $\lambda, \theta$, we define $v_t(\theta)=e^{-\theta t}$ and
\[
x_t(\lambda, \theta)=
\begin{cases}
\frac{1}{1+\theta t} & \text{~if~}\lambda=\theta,\\
\frac{(\lambda-\theta)e^{(\lambda-\theta)t}}{\lambda e^{(\lambda-\theta)t}-\theta} & \text{~else}.
\end{cases}
\]
For $0\leq t\leq T_0$, we define $M_4(t)$ as the $2\times 2$ symmetric matrix such that
\[
M_4(t)=\sum_{i=1}^3l_iF_i\left(x_t(\lambda, \theta), v_t(\theta)\right)l_i^{\mathsf{T}}
\]
and $M_5(t)$ as the $2\times 2$ matrix such that
\[
M_5(t)=\left(\sum_{i=1}^3l_i\nabla^{\mathsf{T}}F_i\right)\left(x_t(\lambda, \theta), v_t(\theta)\right),
\]
where $\nabla=(\frac{\partial}{\partial x}, \frac{\partial}{\partial y})^{\mathsf{T}}$. Then we let $\{Y_t\}_{t\geq 0}$ be the time-inhomogeneous $2-$dimensional O-U process such that
\[
\begin{cases}
&dY_t=M_5(t)Y_tdt+M_4^{\frac{1}{2}}(t)dB_t, \\
&Y_0=0,
\end{cases}
\]
where $\{B_t\}_{t\geq 0}$ is a standard $2-$dimensional Brownian motion. As a result, $Y_t$ follows a Gaussian distribution $\mathbb{N}(\mathbf{0}, \Sigma_t)$ for all $0<t\leq T_0$, where $\Sigma_t$ is a $2\times 2$ positive definite matrix for every $t$. Now we can give our central limit theorem.

\begin{theorem}\label{theorem 2.2 CLT}
As $n$ grows to infinity,
\[
\left(\sqrt{n}\left(\hat{\lambda}_n-\lambda\right), \sqrt{n}\left(\hat{\theta}_n-\theta\right)\right)^{\mathsf{T}}
\]
converges in distribution to $\mathbb{N}\left(\mathbf{0}, M_3\Sigma_{T_0}M_3^{\mathsf{T}}\right)$.
\end{theorem}

\begin{remark}\label{remak 2.1}
We can utilize the following approach introduced in Chapter 11 of \cite{Ethier1986} to approximate $\Sigma_{T_0}$ via a computer. Let $\{\Phi(t)\}_{t\geq 0}$ be the solution to the ODE
\[
\begin{cases}
&\frac{d}{dt}\Phi(t)=-\Phi(t)M_5(t), \\
&\Phi(0)=
\begin{pmatrix}
1 & 0\\
0& 1
\end{pmatrix}
,
\end{cases}
\]
which can be simulated by Euler's method, then, according to Ito's formula,
\[
\Phi(t)Y_t=\int_0^t \Phi(s)M_4^{\frac{1}{2}}(s)dB_s
\]
and hence $\Sigma_{T_0}$ is given by
\[
\int_0^{T_0}\Phi^{-1}(T_0)\Phi(s)M_4(s)\Phi^{\mathsf{T}}(s)\left(\Phi^{-1}(T_0)\right)^{\mathsf{T}}ds.
\]

\end{remark}

Theorem \ref{theorem 2.2 CLT} can be utilized in hypothesis testings of $\lambda$ and $\theta$. For example, let $\lambda_0$ be a known given constant and we discuss the hypothesis testing
\begin{equation}\label{equ hypothesis testing}
H_0: \lambda=\lambda_0\text{\quad}{\rm vs}\text{\quad} H_1: \lambda\neq \lambda_0.
\end{equation}
Since $\left(M_3\Sigma_{T_0}M^{\mathsf{T}}_3\right)(1,1)$ relies on $(\lambda, \theta)$, we write it as $\left(M_3\Sigma_{T_0}M^{\mathsf{T}}_3\right)(1,1, \lambda, \theta)$. Then, by Theorems \ref{theorem 2.1 LLN} and \ref{theorem 2.2 CLT}, under $H_0$,
\[
\frac{\sqrt{n}\left(\hat{\lambda}_n-\lambda_0\right)}{\sqrt{\left(M_3\Sigma_{T_0}M^{\mathsf{T}}_3\right)(1,1, \lambda_0, \hat{\theta}_n)}}
\]
approximately follows standard Normal distribution $\mathbb{N}(0,1)$ for large $n$. Consequently, let
\[
W=\left\{\left|\frac{\sqrt{n}\left(\hat{\lambda}_n-\lambda_0\right)}{\sqrt{\left(M_3\Sigma_{T_0}M^{\mathsf{T}}_3\right)(1,1, \lambda_0, \hat{\theta}_n)}}\right|>1.96\right\},
\]
then $W$ is an approximated reject region at significant level $0.05$.

Our last result is about moderate deviation principles of $\hat{\lambda}_n$ and $\hat{\theta}_n$. Let $\{a_n\}_{n\geq 1}$ be a given positive sequence such that
$\lim_{n\rightarrow+\infty}\frac{a_n}{n}=0$ while $\lim_{n\rightarrow+\infty}\frac{a_n}{\sqrt{n}}=+\infty$ (e.g. $a_n=n^{2/3}$), then we have the following result.

\begin{theorem}\label{theorem 2.3 MDP}
For any $\epsilon>0$, there exists $I_1(\epsilon), I_2(\epsilon)>0$ such that
\[
\lim_{n\rightarrow +\infty}\frac{n}{a_n^2}\log P\left(\left|\frac{n(\hat{\lambda}_n-\lambda)}{a_n}\right|>\epsilon\right)=-I_1(\epsilon)
\]
and
\[
\lim_{n\rightarrow +\infty}\frac{n}{a_n^2}\log P\left(\left|\frac{n(\hat{\theta}_n-\theta)}{a_n}\right|>\epsilon\right)=-I_2(\epsilon).
\]
Furthermore, $I_1(\epsilon)$ and $I_2(\epsilon)$ are given by
\[
I_1(\epsilon)=\frac{\epsilon^2}{2\left(M_3\Sigma_{T_0}M_3^{\mathsf{T}}\right)(1,1,\lambda, \theta)}
\text{~and~}I_2(\epsilon)=\frac{\epsilon^2}{2\left(M_3\Sigma_{T_0}M_3^{\mathsf{T}}\right)(2,2,\lambda, \theta)}.
\]
\end{theorem}

We can utilize Theorem \ref{theorem 2.3 MDP} to give confidence intervals of $\lambda$ and $\theta$. For example, let $a_n=n^\alpha$ with $\alpha\in(\frac{1}{2}, 1)$ and $\epsilon=1$, then, by Theorem \ref{theorem 2.3 MDP},
\[
\left[\hat{\lambda}_n-n^{-(1-\alpha)}, ~\hat{\lambda}_n+n^{-(1-\alpha)}\right]
\]
is a confidence interval of $\lambda$ at confidence level about $1-e^{-I_1(1)n^{2\alpha-1}}$. Note that the above confidence interval has the advantage that the length of the interval grows to $0$ meanwhile the confidence level grows to $1$ exponentially as $n\rightarrow +\infty$. Following are simulation results of the above confidence intervals for $\lambda=2, \theta=1, T_0=1$, $20\leq n\leq 1000$ and $\alpha=\frac{3}{4}, \frac{2}{3}, \frac{3}{5}, \frac{11}{20}$ respectively.

\begin{figure}[H]
\centering
\includegraphics[scale=0.5]{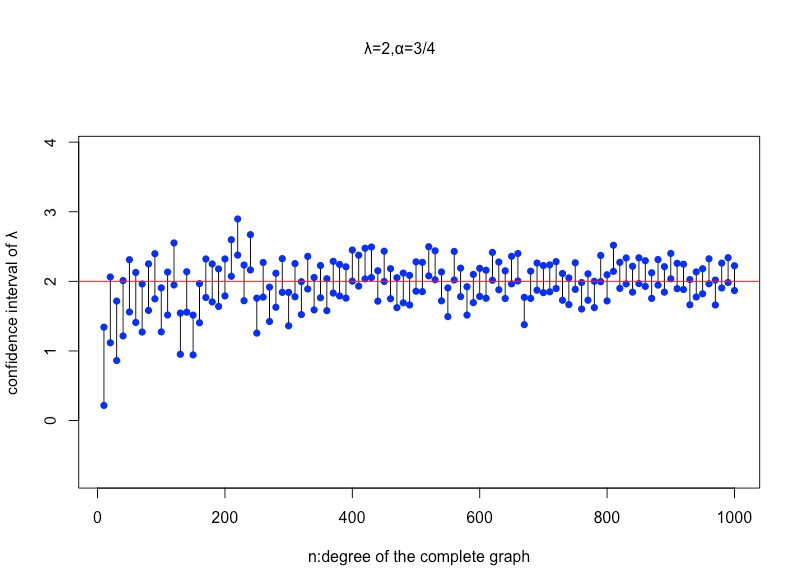}
\caption{$\lambda=2$, $\theta=1$ and $\alpha=3/4$}\label{figure 2.5}
\end{figure}

\begin{figure}[H]
\centering
\includegraphics[scale=0.5]{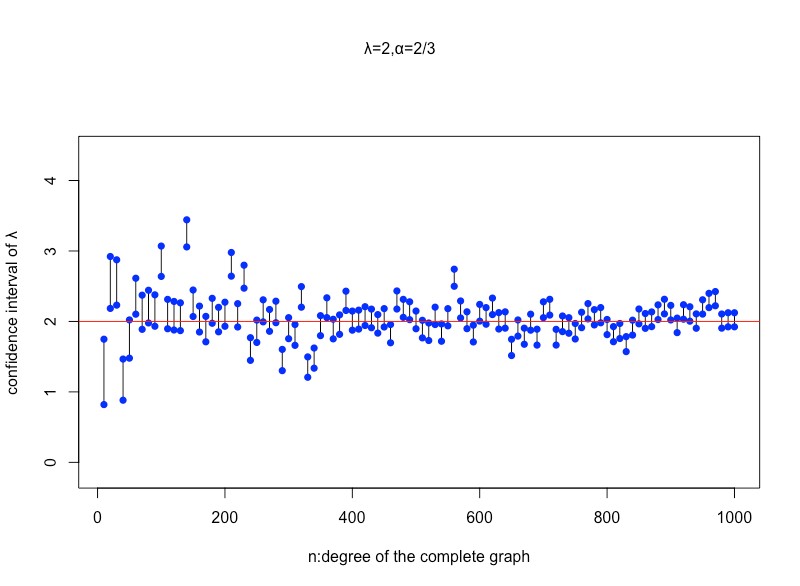}
\caption{$\lambda=2$, $\theta=1$ and $\alpha=2/3$}\label{figure 2.6}
\end{figure}

\begin{figure}[H]
\centering
\includegraphics[scale=0.5]{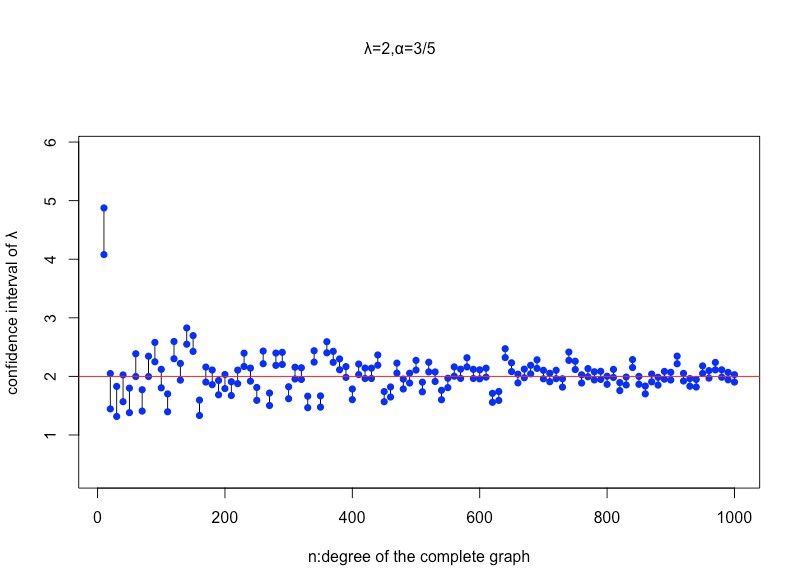}
\caption{$\lambda=2$, $\theta=1$ and $\alpha=3/5$}\label{figure 2.7}
\end{figure}

\begin{figure}[H]
\centering
\includegraphics[scale=0.5]{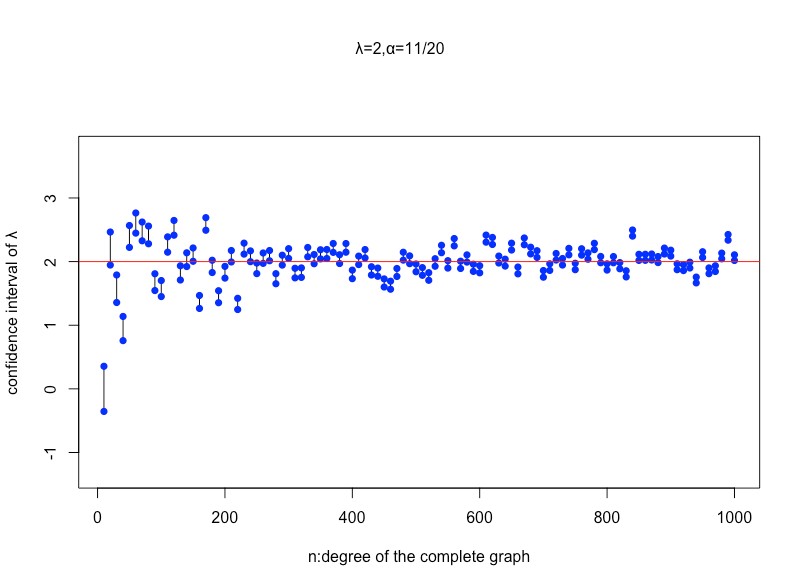}
\caption{$\lambda=2$, $\theta=1$ and $\alpha=11/20$}\label{figure 2.8}
\end{figure}

\section{Outlines of the proofs}\label{section three}

In this section, we give outlines of the proofs of our main theorems. We mainly utilize the theory of density-dependent Markov chains introduced in \cite{Kurtz1978} and Chapter 11 of \cite{Ethier1986} authored by Ethier and Kurtz. As a preparation, we recall the definition of density-dependent Markov processes. For each integer $n\geq 1$, the density-dependent Markov process $\{\xi_t^n\}_{t\geq 0}$ is a continuous-time Markov process with state space $\mathbb{R}^d$ for some $d$. The transition rates functions of $\{\xi_t^n\}_{t\geq 0}$ is given by
\[
\xi_t^n\rightarrow \xi_t^n+l \text{~at rate~}nF_l\left(\frac{\xi_t^n}{n}\right)
\]
for any $l\in \mathcal{D}$, where $\mathcal{D}$ is a given subset of $\mathbb{R}^d$ while $\{F_l\}_{l\in \mathcal{D}}$ are smooth functions from $\mathbb{R}^d$ to $[0, +\infty)$. To illustrate the relation between density-dependent Markov chains and our SIS epidemic model, we define
\[
X_t^n=\left|A_t^n\right|\text{~and~}V_t^n=\left|\left\{x\in C_n:~x\in A_s^n\text{~for all~}0\leq s\leq t\right\}\right|.
\]
Specially, $X_{T_0}^n=X_n$ and $V_{T_0}^n=V_n$, where $X_n, V_n$ are defined as in Section \ref{section two}. Let $\zeta_t^n=\left(X_t^n, V_t^n\right)^{\mathsf{T}}$, then, according to the transition rates functions of $\{A_t^n\}_{t\geq 0}$, it is easy to check that $\{\zeta_t^n\}_{t\geq 0}$ is a version of density-dependent Markov chain with $d=2$,
\[
\mathcal{D}=\left\{l_1, l_2, l_3\right\}
\]
and $F_{l_i}=F_i$ for $i=1,2,3$, where $l_1, l_2, l_3$ and $F_1, F_2, F_3$ are defined as in Section \ref{section two}.

We recall the law of large numbers and cental limit theorem of density-dependent Markov chains introduced in \cite{Kurtz1978} by Kurtz. Note that vectors appear in following propositions are all column vectors.

\begin{proposition}\label{Kurtz LLN} (Kurtz, 1978)
If $\xi_0^n=n\xi_0$ for all $n\geq 1$, then $\frac{\xi_t^n}{n}$ converges in probability to the solution to the ODE
\begin{equation}\label{equ 3.1 ODE}
\begin{cases}
&\frac{d}{dt}x_t=\sum_{l\in \mathcal{D}}lF_l(x_t),\\
&x_0=\xi_0
\end{cases}
\end{equation}
as $n$ grows to infinity.
\end{proposition}

\begin{proposition}\label{Kurtz CLT}(Kurtz, 1978)
If $\frac{\xi_0^n-n\xi_0}{\sqrt{n}}=w_0$ for all $n\geq 1$, then $\frac{\xi_t^n-nx_t}{\sqrt{n}}$ converges in distribution to the time-inhomogeneous O-U process
\[
\begin{cases}
&dW_t=M_6(t)W_tdt+M_7^{\frac{1}{2}}(t)dB_t,\\
&W_0=w_0
\end{cases}
\]
as $n$ grows to infinity, where $x_t$ is the solution to Equation \eqref{equ 3.1 ODE},
\[
M_6(t)=\sum_{l\in \mathcal{D}}\left(l\nabla^{\mathsf{T}}F_l(x_t)\right), \text{~}M_7(t)=\sum_{l\in\mathcal{D}}lF_l\left(x_t\right)l^{\mathsf{T}}
\]
and $\{B_t\}_{t\geq 0}$ are $d$-dimensional standard Brownian motions.
\end{proposition}

Now we can give proofs of our first and second main results.

\proof[Proof of Theorem \ref{theorem 2.1 LLN}]

By Proposition \ref{Kurtz LLN}, $\left(\frac{X_n}{n}, \frac{V_n}{n}\right)^{\mathsf{T}}$ converges in probability to $\left(x_{_{T_0}}, v_{_{T_0}}\right)^{\mathsf{T}}$, where
$(x_t, v_t)^{\mathsf{T}}$ is the solution to
\begin{equation}\label{equ 3.3 ode}
\begin{cases}
&\frac{d}{dt}x_t=-\theta x_t+\lambda x_t(1-x_t),\\
&\frac{d}{dt}v_t=-\theta v_t,\\
&\left(x_0, v_0\right)^{\mathsf{T}}=(1,1)^{\mathsf{T}}.
\end{cases}
\end{equation}
By directly solving the above Equation, $x_{_{T_0}}=K(\theta, \lambda)$ while $v_{_{T_0}}=e^{-\theta T_0}$. Therefore,
\[
\lambda=H(\theta, x_{_{T_0}})\text{~and~}\theta=-\frac{1}{T_0}\log v_{_{T_0}},
\]
Theorem \ref{theorem 2.1 LLN} follows from which directly since $H(x,y)$ and $\log x$ are continuous functions.

\qed

\proof[Proof of Theorem \ref{theorem 2.2 CLT}]

For simplicity, we use $o_p(1)$ to denote a random variable $\varepsilon_n$ when $\varepsilon_n$ converges in probability to $0$ as $n\rightarrow +\infty$. By Theorem \ref{theorem 2.1 LLN} and Lagrange's mean value theorem,
\begin{align*}
\left(\hat{\lambda}_n-\lambda, \hat{\theta}_n-\theta\right)^{\mathsf{T}}&=\left(H\left(\hat{\theta}_n, \frac{X_n}{n}\right)-H(\theta, x_{_{T_0}}), \hat{\theta}_n-\theta\right)^{\mathsf{T}} \\
&=\left(\begin{pmatrix}
\frac{\partial H}{\partial z}\left(\theta, x_{_{T_0}}\right) & \frac{\partial H}{\partial x}\left(\theta, x_{_{T_0}}\right)\\
0 & 1
\end{pmatrix}+o_p(1)\right)\begin{pmatrix} \frac{X_n}{n}-x_{_{T_0}}\\ \hat{\theta}_n-\theta\end{pmatrix}\\
&=\left(M_1+o_p(1)\right)\left(\frac{X_n}{n}-x_{_{T_0}}, \hat{\theta}_n-\theta\right)^{\mathsf{T}},
\end{align*}
since $x_{_{T_0}}=K(\theta, \lambda)$. According to a similar analysis,
\begin{align*}
\left(\frac{X_n}{n}-x_{_{T_0}}, \hat{\theta}_n-\theta\right)^{\mathsf{T}}&=\left(\frac{X_n}{n}-x_{_{T_0}}, -\frac{1}{T_0}\log \frac{V_n}{n}-\left(-\frac{1}{T_0}\log v_{_{T_0}}\right)\right)^{\mathsf{T}}\\
&=\left(\begin{pmatrix}
1 & 0\\
0 & -\frac{1}{T_0}\frac{1}{v_{_{T_0}}}
\end{pmatrix}+o_p(1)\right)\begin{pmatrix} \frac{X_n}{n}-x_{_{T_0}}\\ \frac{V_n}{n}-v_{_{T_0}}\end{pmatrix}\\
&=\left(M_2+o_p(1)\right)\left(\frac{X_n}{n}-x_{_{T_0}}, \frac{V_n}{n}-v_{_{T_0}}\right)^{\mathsf{T}}.
\end{align*}
Consequently,
\begin{equation}\label{equ 3.2}
\left(\hat{\lambda}_n-\lambda, \hat{\theta}_n-\theta\right)^{\mathsf{T}}=\left(M_3+o_p(1)\right)\left(\frac{X_n}{n}-x_{_{T_0}}, \frac{V_n}{n}-v_{_{T_0}}\right)^{\mathsf{T}}.
\end{equation}
By Proposition \ref{Kurtz CLT}, $\sqrt{n}\left(\frac{X_n}{n}-x_{_{T_0}}, \frac{V_n}{n}-v_{_{T_0}}\right)^{\mathsf{T}}$ converges in distribution to $Y_{T_0}$ as $n\rightarrow+\infty$, where $Y_{T_0}$ is defined as in Section \ref{section two}. That is to say,
\[
\sqrt{n}\left(\frac{X_n}{n}-x_{_{T_0}}, \frac{V_n}{n}-v_{_{T_0}}\right)^{\mathsf{T}}
\]
converges in distribution to $\mathbb{N}\left(\mathbf{0},\Sigma_{T_0}\right)$ as $n\rightarrow +\infty$, Theorem \ref{theorem 2.2 CLT} follows from which and Equation \eqref{equ 3.2} directly.

\qed

Based on Theorem \ref{theorem 2.2 CLT}, readers not familiar with theories of moderate deviations could intuitively understand Theorem \ref{theorem 2.3 MDP} in the following way. Theorem \ref{theorem 2.2 CLT} can be roughly written as
\[
P\left(\sqrt{n}\left(\hat{\lambda}_n-\lambda\right)=dx\right)\approx \exp\left\{-\frac{x^2}{2\left(M_3\Sigma_{T_0}M_3^{\mathsf{T}}\right)(1,1,\lambda,\theta)}\right\}dx.
\]
Then,
\begin{align*}
P\left(\frac{n\left(\hat{\lambda}_n-\lambda\right)}{a_n}=dx\right)&=P\left(\sqrt{n}\left(\hat{\lambda}_n-\lambda\right)=\frac{a_n}{\sqrt{n}}dx\right)\\
&\approx\exp\left\{-\frac{a_n^2x^2}{2n\left(M_3\Sigma_{T_0}M_3^{\mathsf{T}}\right)(1,1,\lambda,\theta)}\right\}dx,
\end{align*}
i.e.,
\[
\lim_{n\rightarrow+\infty}\frac{n}{a_n^2}\log P\left(\frac{n\left(\hat{\lambda}_n-\lambda\right)}{a_n}=dx\right)=-I_1(x).
\]
The rigorous proof of Theorem \ref{theorem 2.3 MDP} is given in the appendix, where a moderate deviation principle for density-dependent Markov chains given in \cite{Xue2018} is utilized. Readers who are convinced by the above intuitive explanation and not interested in too many mathematical details could just skip this proof.

\appendix
\section{Proof of Theorem \ref{theorem 2.3 MDP}}

\proof[Proof of Theorem \ref{theorem 2.3 MDP}]

 We denote by $\mathcal{S}$ the set of functions from $[0, T_0]$ to $\mathbb{R}^2$ which are right continuous, have left-hand limits and starts at $(0,0)^{\mathsf{T}}$, i.e., the set of c\`{a}dl\`{a}g functions $f$ with $f(0)=(0,0)^{\mathsf{T}}$. Let $\zeta_t^n=\left(X_t^n, V_t^n\right)^{\mathsf{T}}$ be defined as in Section \ref{section three} while $(x_t, v_t)^{\mathsf{T}}$ be the solution to Equation \eqref{equ 3.3 ode}, then, by Theorem 2.1 of \cite{Xue2018}, the path $\vartheta_n:=\left\{\frac{n}{a_n}\left(\frac{\zeta_t^n}{n}-(x_t, v_t)^{\mathsf{T}}\right)\right\}_{0\leq t\leq T_0}$ follows moderate deviation principle with rate function $J_1(\cdot)$ given by
 \[
 J_1(f)=
 \begin{cases}
 &\frac{1}{2}\int_0^{T_0}\left(f^\prime_t-M_5(t)f_t\right)^{\mathsf{T}}M_4^{-1}(t)\left(f^\prime_t-M_5(t)f_t\right)dt\\
 & \text{~\quad\quad\quad if $f$ is absolutely continuous},\\
 & +\infty  \text{~else}.
 \end{cases}
 \]
 That is to say,
 \[
 \limsup_{n\rightarrow+\infty}\frac{a_n^2}{n}\log P\left(\vartheta_n\in C\right)\leq -\inf_{f\in C}J_1(f)
 \]
for any closed set $C\subseteq \mathcal{S}$ while
\[
\liminf_{n\rightarrow+\infty}\frac{a_n^2}{n}\log P\left(\vartheta_n\in C\right)\geq -\inf_{f\in O}J_1(f)
\]
for any open set $O\subseteq \mathcal{S}$. Then, according to the contraction principle (see Section 4.2 of \cite{Dembo1997} authored by Dembo and Zeitouni),
$\frac{n}{a_n}\left(\frac{X^n_{T_0}}{n}-x_{_{T_0}}, \frac{V^n_{T_0}}{n}-v_{_{T_0}}\right)^{\mathsf{T}}$ follows moderate deviation principle with rate function $J_2(\cdot)$ given by
\[
J_2(x)=\inf_{f\in \mathcal{S}, f_{_{T_0}}=x}J_1(f)
\]
for any $x=(x_1, x_2)^{\mathsf{T}}\in \mathbb{R}^2$. We claim that
\begin{equation}\label{equ A.1}
J_2(x)=\frac{x^{\mathsf{T}}\Sigma^{-1}_{T_0}x}{2}
\end{equation}
for any $x\in \mathbb{R}^2$. Equation \eqref{equ A.1} holds according to an utilization of Cauchy-Schwartz inequality, the detail of which we put at the end of this appendix.

According to the analysis given in the proof of Theorem \ref{theorem 2.2 CLT},
\[
\frac{n}{a_n}\left(\hat{\lambda}_n-\lambda, \hat{\theta}_n-\theta\right)^{\mathsf{T}}=\frac{n}{a_n}\left(M_3+\epsilon_n\right)\left(\frac{X^n_{T_0}}{n}-x_{_{T_0}}, \frac{V^n_{T_0}}{n}-v_{_{T_0}}\right)^{\mathsf{T}},
\]
where $\epsilon_n=o_p(1)$. According to large deviation principles of epidemic models established in \cite{Pardoux2017} by Pardoux and Samegni-Kepgnou, $\epsilon_n$ follows a large deviation principle with a rate function $I_3(\cdot)$, i.e., for any $\epsilon>0$, there exists $I_3(\epsilon)>0$ such that
\[
\limsup_{n\rightarrow+\infty}\frac{1}{n}\log P\left(\sum_{i=1}^2\sum_{j=1}^2\left|\epsilon_n(i,j)\right|\geq \epsilon\right)\leq -I_3(\epsilon).
\]
Consequently, since $\frac{a_n^2}{n}=o(n)$,
\[
\lim_{n\rightarrow+\infty}\frac{n}{a_n^2}\log P\left(\sum_{i=1}^2\sum_{j=1}^2\left|\epsilon_n(i,j)\right|\geq \epsilon\right)=-\infty
\]
and hence $\frac{n}{a_n}\left(\hat{\lambda}_n-\lambda, \hat{\theta}_n-\theta\right)^{\mathsf{T}}$ and $\frac{n}{a_n}M_3\left(\frac{X^n_{T_0}}{n}-x_{_{T_0}}, \frac{V^n_{T_0}}{n}-v_{_{T_0}}\right)^{\mathsf{T}}$ follows the same moderate deviation principle. As a result, by the contraction principle, Theorem \ref{theorem 2.3 MDP} holds with
\[
I_1(\epsilon)=\inf_{M_3(1,1)x_1+M_3(1,2)x_2=\epsilon}J_2(x)=\inf_{M_3(1,1)x_1+M_3(1,2)x_2=\epsilon}\frac{x^{\mathsf{T}}\Sigma^{-1}_{T_0}x}{2}
\]
while
\[
I_2(\epsilon)=\inf_{M_3(2,1)x_1+M_3(2,2)x_2=\epsilon}J_2(x)=\inf_{M_3(2,1)x_1+M_3(3,2)x_2=\epsilon}\frac{x^{\mathsf{T}}\Sigma^{-1}_{T_0}x}{2}.
\]
Let $\eta_1=\left(M_3(1,1), M_3(1,2)\right)^{\mathsf{T}}$ and $y=\frac{\epsilon\Sigma_{T_0}\eta_1}{\eta_1^T\Sigma_{T_0}\eta_1}$, then $M_3(1,1)y_1+M_3(1,2)y_2=\epsilon$ and $\left(x-y\right)^{\mathsf{T}}\Sigma^{-1}_{T_0}y=0$ for any $x\in \mathbb{R}^2$ satisfying $\eta_1^T x=\epsilon$. Therefore,
\[
\frac{x^{\mathsf{T}}\Sigma^{-1}_{T_0}x}{2}=\frac{y^{\mathsf{T}}\Sigma^{-1}_{T_0}y}{2}+\frac{(x-y)^{\mathsf{T}}\Sigma^{-1}_{T_0}(x-y)}{2}
\]
for any $x$ satisfying $\eta_1^{\mathsf{T}}x=\epsilon$. Then, since $\Sigma_{T_0}$ is positive definite,
\begin{align*}
I_1(\epsilon)&=\inf_{\eta_1^{\mathsf{T}}x=\epsilon}\frac{x^{\mathsf{T}}\Sigma^{-1}_{T_0}x}{2}=\frac{y^{\mathsf{T}}\Sigma^{-1}_{T_0}y}{2}
=\frac{\epsilon^2}{2\left(\eta_1^{\mathsf{T}}\Sigma_{T_0}\eta_1\right)}=\frac{\epsilon^2}{2\left(M_3\Sigma_{T_0}M_3^{\mathsf{T}}\right)(1,1,\lambda, \theta)}.
\end{align*}
According to a similar analysis,
\[
I_2(\epsilon)=\frac{\epsilon^2}{2\left(M_3\Sigma_{T_0}M_3^{\mathsf{T}}\right)(2,2,\lambda, \theta)}
\]
and hence the proof is complete.

\qed

At last, we only need to prove Equation \eqref{equ A.1}.

\proof[Proof of Equation \eqref{equ A.1}]

Let $\Phi(t)$ be defined as in Remark \ref{remak 2.1}, then, the ODE
\[
\begin{cases}
&f^\prime_t-M_5(t)f_t=h(t), \\
&f(0)=(0,0)^{\mathsf{T}}
\end{cases}
\]
has the unique solution
\begin{equation}\label{equ A.2}
f_t=\int_0^t\Phi^{-1}(t)\Phi(s)h(s),
\end{equation}
since
\[
\Phi(t)f^{\prime}_t-\Phi(t)M_5(t)f_t=\left(\Phi(t)f_t\right)^{\prime}.
\]
For later use, we need choose a $h(t)$ with form $h(t)=M_4(t)\Phi^{\mathsf{T}}(t)b$ for some $b\in \mathbb{R}^2$ to make $f$ given by \eqref{equ A.2} satisfy $f_{_{T_0}}=x$.  By direct calculation, we let
\[
b=\left(\int_0^{T_0}\Phi^{-1}(T_0)\Phi(s)M_4(s)\Phi^{-1}(s)ds\right)^{-1}x=\left(\Phi^{-1}(T_0)\right)^{\mathsf{T}}\Sigma^{-1}_{T_0}x
\]
and then
\[
h(t)=M_4(t)\Phi^{\mathsf{T}}(t)\left(\Phi^{-1}(T_0)\right)^{\mathsf{T}}\Sigma^{-1}_{T_0}x,
\]
since $\Sigma_{T_0}=\int_0^{T_0}\Phi^{-1}(T_0)\Phi(s)M_4(s)\Phi^{\mathsf{T}}(s)\left(\Phi^{-1}(T_0)\right)^{\mathsf{T}}ds$ as we have shown in Remark \ref{remak 2.1}.  For this $h(t)$, let $f$ be defined as in \eqref{equ A.2}, then
\begin{align*}
J_1(f)&=\frac{\int_0^{T_0}h^{\mathsf{T}}(t)M_4^{-1}(t)h(t)dt}{2} \\
&=\frac{x^{\mathsf{T}}\Sigma^{-1}_{T_0}\left(\int_0^{T_0}\Phi^{-1}(T_0)\Phi(t)M_4(t)\Phi^{\mathsf{T}}(t)\left(\Phi^{-1}(T_0)\right)^{\mathsf{T}}dt\right)\Sigma^{-1}_{T_0}x}{2}\\
&=\frac{x^{\mathsf{T}}\Sigma^{-1}_{T_0}\Sigma_{T_0}\Sigma^{-1}_{T_0}x}{2}=\frac{x^{\mathsf{T}}\Sigma^{-1}_{T_0}x}{2}.
\end{align*}
As a result,
\begin{equation}\label{equ A.3}
J_2(x)\leq \frac{x^{\mathsf{T}}\Sigma^{-1}_{T_0}x}{2}.
\end{equation}
On the other hand, for any absolutely continuous $g\in \mathcal{S}$ and any $k\in \mathcal{S}$,
\begin{align*}
&J_1(g)\int_0^{T_0}k_t^{\mathsf{T}}k_t dt\\
&=\frac{1}{2}\int_0^{T_0}k_t^{\mathsf{T}}k_t dt\int_0^{T_0}\left(M^{-\frac{1}{2}}_4(t)\left(g^\prime_t-M_5(t)g_t\right)\right)^{\mathsf{T}}\left(M^{-\frac{1}{2}}_4(t)\left(g^\prime_t-M_5(t)g_t\right)\right) dt\\
&\geq \frac{1}{2}\left(\int_0^{T_0} k_t^{\mathsf{T}}M_4^{-\frac{1}{2}}(t)\left(g^\prime_t-M_5(t)g_t\right)dt\right)^2
\end{align*}
according to Cauchy-Schwartz inequality. We choose $k_t=M_4^{\frac{1}{2}}(t)\Phi^{\mathsf{T}}(t)b$, then, for $g$ satisfying $g_{_{T_0}}=x$,
\[
\left(\int_0^{T_0} k_t^{\mathsf{T}}M_4^{-\frac{1}{2}}(t)\left(g^\prime_t-M_5(t)g_t\right)dt\right)^2=
\left(b^{\mathsf{T}}\int_0^{T_0}\left(\Phi(t)g_t\right)^{\prime}dt\right)^2=\left(b^{\mathsf{T}}\Phi(T_0)x\right)^2.
\]
As a result,
\[
J_2(x)\geq \frac{\left(b^{\mathsf{T}}\Phi(T_0)x\right)^2}{2b^{\mathsf{T}}\left(\int_0^{T_0}\Phi(t)M_4(t)\Phi^{\mathsf{T}}(t)dt\right)b}.
\]
According to the definition of $b$,
\[
\frac{\left(b^{\mathsf{T}}\Phi(T_0)x\right)^2}{2b^{\mathsf{T}}\left(\int_0^{T_0}\Phi(t)M_4(t)\Phi^{\mathsf{T}}(t)dt\right)b}
=\frac{\left(x^{\mathsf{T}}\Sigma^{-1}_{T_0}x\right)^2}{2x^{\mathsf{T}}\Sigma^{-1}_{T_0}x}=\frac{x^{\mathsf{T}}\Sigma^{-1}_{T_0}x}{2}.
\]
Therefore, $J_2(x)\geq \frac{x^{\mathsf{T}}\Sigma^{-1}_{T_0}x}{2}$, Equation \eqref{equ A.1} follows directly from which and \eqref{equ A.3}.

\qed

\quad

\textbf{Acknowledgments.} The authors are grateful to the financial support from the National Natural Science Foundation of China with
grant numbers 11501542.

{}

\begin{thebibliography}{}
\bibitem{Becker1977}Becker, N. (1977).  Estimation for discrete time branching processes with application to epidemics.  \emph{Biometrics} \textbf{33}, 515-522.
\bibitem{Dembo1997}Dembo, A. and Zeitouni, O. (1997). \emph{Large Deviations: Techniques and Applications.} Springer, Berlin.
\bibitem{Ethier1986}Ethier, N. and Kurtz, T. (1986). \emph{Markov Processes: Characterization and Convergence.} John Wiley and Sons, Hoboken, NJ, USA.
\bibitem{Fierro2015}Fierro, R., Leiva, V. and Balakrishnan, N. (2015). Statistical inference on a stochastic epidemic model. \emph{Communications in Statistics. Simulation and Computation} \textbf{44}, 2297-2314.
\bibitem{Guy2015} Guy, R., Lar\'{e}do, C. and Vergu, E. (2015). Approximation of epidemic models by diffusion processes and their statistical inference.  \emph{Journal of Mathematical Biology} \textbf{70}, 621-646.
\bibitem{Hadeler2011}Hadeler, K. P. (2011). Parameter estimation in epidemic models: simplified formulas. \emph{Canadian Applied Mathematics Quarterly} \textbf{19}, 343-356.
\bibitem{Kurtz1978}Kurtz, T. (1978). Strong approximation theorems for density dependent Markov chains. \emph{Stochastic Processes and their Applications} \textbf{6}, 223-240.
\bibitem{Lekone2006}Lekone, P., Finkenst\"{a}dt, B. (2006). Statistical inference in a stochastic epidemic SEIR model with control intervention: Ebola as a case study. \emph{Biometrics} \textbf{62}, 1170-1177.
\bibitem{Lig1985}Liggett, T. M. (1985). {\it Interacting Particle Systems.} Springer, New York.
\bibitem{Lin2013}Lindenstrand, D. and Svensson, {\AA}. (2013). Estimation of the Malthusian parameter in an stochastic epidemic model using martingale methods. \emph{Mathematical Biosciences} \textbf{246}, 272-279.
\bibitem{Pan2014}Pan, J., Gray, A., Greenhalgh, D. and Mao, X. (2014). Parameter estimation for the stochastic SIS epidemic model. \emph{Statistical Inference for Stochastic Processes} \textbf{17}, 75-98.
\bibitem{Pardoux2017}Pardoux, E. and Samegni-Kepgnou, B. (2017). Large deviation principle for epidemic models. \emph{Journal of Applied Probability} \textbf{54}, 905-920.
\bibitem{Xue2018}Xue, XF. (2018). Moderate deviations of density-dependent Markov chains. Arxiv: 1908.03762.
\bibitem{Yip1998} Yip, P. and Chen, Q. (1998). Statistical inference for a multitype epidemic model. \emph{Journal of Statistical Planning and Inference} \textbf{71}, 229-244.
\end{thebibliography}
\end{document}